\documentclass[11pt]{article}
\usepackage{hyperref}
\usepackage[latin1]{inputenc}
\usepackage{amsmath}
\usepackage{amsfonts}
\usepackage{amsthm}
\usepackage{amssymb}

\usepackage{hhline}
\usepackage{array} 
\newcolumntype{M}[1]{>{\centering\arraybackslash}m{#1}}
\newcolumntype{N}{@{}m{0pt}@{}}

\usepackage{color}
\usepackage[left=4.00cm, right=3.00cm, top=4.00cm, bottom=3.00cm]{geometry}
\usepackage{graphicx}
\usepackage{pdfsync}
\usepackage{tikz}
\usepackage{undertilde}
\usepackage[active]{srcltx}
\usepackage{pgf}
\usepackage{mathtools}
\usetikzlibrary{shadows}
\usepackage{accents}
\usepackage{todonotes}
\usepackage{bm}
\newtheorem{defn}{Definition}[section]
\newtheorem{thm}[defn]{Theorem}

\newtheorem{lem}[defn]{Lemma}
\newtheorem{cor}[defn]{Corollary}

\newtheorem{remark}{Remark}
\newtheorem{assumption}[defn]{Assumption}

\newfont{\klein}{cmcsc10 scaled 900}
\newfont{\mbs}{cmmib7}
\newfont{\mbbs}{cmmib5}
\def\DS{\mbox{{\boldmath $S$}}}
\def\Sm{{\mbs S}}
\def\Ss{{\mbbs S}}
\def\state{{\mathchoice{\DS}{\DS}{\Sm}{\Ss}}}
\def\DZZ{{\hbox{\bf Z}}}
\def\SZZ{{\hbox{\bf\scriptsize Z}}}
\def\SSZZ{{\hbox{\bf\tiny Z}}}
\def\intgr{{\mathchoice\DZZ\DZZ\SZZ\SSZZ}}
\newcommand{\stless}{\stackrel{st}{\preceq}}
\def\DRe{{\hbox{\bf R}}}
\def\SRe{{\hbox{\scriptsize R}}}
\def\SSRe{{\hbox{\tiny R}}}
\def\Re{{\mathchoice\DRe\DRe\SRe\SSRe}}

\newcommand{\inp}{\cdot}
\newcommand\gn[2][{}]{\prescript{}{#1}q_{#2}}
\newcommand\Gn[1][{}]{\prescript{}{#1}\!Q}
\newcommand\PrT[2][{}]{\prescript{}{#1}\!P^{#2}}
\newcommand\Prt[3][{}]{\,\null_{\raisebox{-.05truecm}{$\scriptstyle #1$}}p_{#2}^{#3}}

\newcommand\tT[2][{}]{\,\null_{#1}\!\tau_{#2}}

 \newcommand{\tcT}[1][{}]{%
   \def\cmdDefaultA{\,{\null^{#1}}}%
   \setbox1=\hbox{$\scriptstyle{#1}\hspace*{.1em}$}
   \cmdNext
 }
\newcommand\tuD{{\underline{D}}} 
 \newcommand{\cmdNext}[1][{}]{%
\setbox2=\hbox{$\scriptstyle{#1}$}\ifdim \wd2<\wd1\setbox2=\hbox to \wd1{\hfill$\scriptstyle{#1}$\hfill}\fi
\ifdim\wd2>\wd1\hspace*{.6\wd2}\fi
{\cmdDefaultA}_{\llap{$\box2$}}\!\!{\sf T}
}
\newcommand\Exp[1][{}]{{\sf E}_{\raisebox{-1pt}{$\scriptstyle #1$}}}
\newcommand\Prob[1][{$\,$}]{{\sf P}_{\raisebox{-2pt}{%
\rlap{$\scriptstyle \!\!#1$}\phantom{$\scriptstyle\!\!\!\!\!#1$}
}}}

\newcommand{\elabel}[1]{\label{#1}%
\vbox to 0pt{\vss\hbox{\hspace{-3cm}\tiny (e:#1)}}}

\newcommand{\Lo}{\widetilde{L}}
\newcommand{\Lu}{\utilde{L}}
\def\ident{\hbox{\sf I}}
\def\elabel#1{\label{#1}}
\def\llabel#1{\label{#1}}
\def\eq #1/{(\ref{#1})}
\def\beq{\begin{equation}}
\def\eeq{\end{equation}}
\catcode`\@=11
\@addtoreset{equation}{section}
\catcode`\@=12

\long\def\notes#1{\ifinner
             {\tiny #1}
             \else
              \marginpar{\protect\tiny #1}
              \fi}

\title{Level product form QSF processes\\
and an analysis of queues with Coxian inter-arrival distribution}
\author{D. Ertiningsih$^{(1)(2)}$ \quad M.N. Katehakis$^{(3)}$ \quad L.C. Smit$^{(1)(3)}$ \quad F.M. Spieksma$^{(1)}$\\
{\small Leiden University, The Netherlands $^{(1)}$ \quad Universitas Gadjah Mada, Indonesia $^{(2)}$}\\ {\small Rutgers Business School, USA $^{(3)}$ }}
\begin{document}
\maketitle
\date
\begin{abstract}
In this paper we study a class of Quasi-Skipfree (QSF) processes where the transition rate submatrices in the skipfree direction have a column times row structure.
Under homogeneity and irreducibility assumptions we show that the stationary distributions of these  processes have a product form as a function of the level.
For an application, we will discuss the ${\it Cox(k)}/M^Y/1$-queue, that can be modelled as a QSF process on a two-dimensional state space. In addition we study the properties of the stationary distribution and 
derive    
monotonicity   of the mean number of the customers in the queue, their mean sojourn time and the variance as a function of $k$ for fixed mean arrival rate.

\vspace*{\baselineskip}\noindent

\vspace*{\baselineskip}\noindent
{\klein AMS Subject Classification} 
\end{abstract}

\section{Introduction}

A Quasi-Skipfree (QSF) process is a continuous time Markov process $X=\{X_t\}_{t\geq 0}$, on a two-dimensional state space $\state=\{(m,i)|\,m\in\intgr,\, i\in\{0,\ldots,\ell_m\}\}$, where $m$ denotes the `level' of the state and $i$ denotes the `phase' within the level. Additionally, the jump rates are not allowed to cross more than one level in one direction, i.e. either the {\em downward} direction or {\em upward} direction. 
This framework is the natural extension of the embedded GI/M/1 and M/G/1-queues, and it has interesting structural properties in common with these processes.

Neuts \cite{neuts1981matrix} investigates the embedded GI/M/1-queue as a skip-free to the right process. The matrix-geometric method for computing the stationary distribution of  homogeneous Quasi-Birth-Death (QBD) processes has been applied in   \cite{neuts1981matrix} and    \cite{latouche1987introduction}. A homogeneous QBD process is a special case of  homogeneous QSF processes, where the transition probabilities are not allowed to cross more than one level in both directions. In his book, Neuts discusses some examples of  homogeneous QBD processes, e.g., the $M/PH/1$-queue and the $PH/M/c$-queue. He shows that the stationary distribution of these queues can be expressed in terms of a rate matrix.

Recently, Katehakis and Smit \cite{katehakis2012successive} have introduced a new procedure to compute the invariant measure for a class of Markov chains. This procedure is called  the successive lumping method. Further, an explicit solution and bounds for the steady state probabilities for the class of ergodic QSF processes that possess the successive lumpability property have been derived in the paper of Katehakis, Smit, and Spieksma \cite{kat-smit-sp-I-13}; an analysis  of the applicability requirements and numerical complexity of the 
successive lumping method is given in  
\cite{kat-smit-sp16}. 
Ramaswami and Latouche \cite{ramaswami1986general} discuss QBD processes with a special structure, namely where the upward or the downward transitions rates form a row times column matrix.
In the latter case the rate matrix can be computed explicitly. In the former case they show that the stationary distribution has a level product form,
as is the case with the embedded $GI/M/1$-queue.
Regarding rate matrix analysis, Ramaswami and Lucantoni, \cite{ramaswami1985algorithms}, exploit the structural properties of the $G|PH|1$-queue, and extend the numerical feasibility of the matrix geometric approach to a wide set of problems by developing efficient  schemes to compute the rate matrix.

One of the main assumptions of \cite{ramaswami1986general} is that the transition rates are bounded as a function of the states.
In the particular case of a Quasi Birth and Death (QBD) process, where the upward transition rates form a matrix with one non-zero row, the rate matrix cannot be computed explicitly. However, we will show how the stationary distribution of a higher level can be explicitly  expressed in terms of the
stationary distribution of lower levels, under an additional invertibility condition. Moreover, we will show that for this specific type of non-homogenous QSF processes it is possible to derive a level product form solution.

We further study  a phase-type inter-arrival, batch service queue, denoted $PH/M^Y/1$-queue, as a particular example of a QSF process, where the upward transition rates form a matrix with one non-zero row.
We will consider a special case of a phase type distribution, a Coxian distribution (cf.~\cite{bhulai2006value}). Herein customers go through a maximum of $k$ exponentially distributed phases, and after each phase the customer can  enter the system. We will denote a Coxian distribution with $k$ phases by $Cox(k)$. Since the Coxian distributions are dense in the set of nonnegative distribution functions, we can approximate those by a Coxian distribution by using for example an expectation-maximization  algorithm (see \cite{asmussen1996fitting}). Therefore the results presented in this paper can be useful for various queues with different inter-arrival distributions.

We show that the stationary distribution of the phases within a level has a product form as well, if the phase-type inter-arrival distribution is a  Coxian distribution. 
In \cite{kat-smit-sp-I-13}  the successive lumpable structure is specified for QSF processes, but we did not investigate the level product form that is derived in this paper for QSF problems with an unbounded number of levels to the left. 
The parameter of the product form is the solution to a fixpoint  equation, that can be numerically approximated efficiently.

We conjecture that allowing for $c$ servers will yield a product form solution based on $c$ factors, similarly to the results presented by Adan et al
\cite[Theorem 4.1]{adan1996analyzing}. Actually, we can also handle the $Cox(k)/E_{l}/1$ by taking the expected workload instead of the number of customers as the batch service distribution.

This paper is organised as follows. In  Section~\ref{s2}  we show that the stationary distribution of the levels has a product form when choosing the downward matrix $D$ as a multiplication of  a column vector $c$ and a row vector $r$. In Section~\ref{s3} we assume that there exists an exit state per level. This is equivalent to the existence of an entrance state per level for a revised level partition of the state space that keeps the QSF property intact. This means that the process is successively lumpable with respect to the new partition. This is used to derive an expression for the rate matrix with respect to the original partition in Section~\ref{subs3}. This derivation is justified by the invertibility of the generator of a transient Markov chain in the Appendix.
Unfortunately, no explicit expression for the rate matrix can be derived. However, we can explicitly determine a matrix that expresses the stationary distribution of a lower level in terms of the stationary distribution of the higher level.

In Section~\ref{s4} we analyze the non-homogenous $Cox(k)/M^{Y}/1$-queue, and find the parameter of the product form, as well as ergodicity properties.
In the remainder of this section we specialise this queue to a queue with homogenous rates and an infinite number of phases.  
In Section~\ref{s5} we derive  monotonicity properties of the stationary mean number of customers in the queue  and their mean sojourn time for fixed mean inter-arrival times. To conclude, we show that the stationary distribution of the $E_k$/M/1-queue (a special case of the $Cox(k)/M^{Y}/1$ where the number of phases is fixed and the batch size of service is $1$) converges monotonically to the stationary distribution of the D/M/1-queue.

\section{The model and basic properties}
\llabel{s2}

In this section we introduce the  notation and derive two initial  properties  
 for the stationary distribution of  homogeneous QSF processes in Lemma \ref{inverse} and Theorem \ref{stat}.
Note that by the homogeneity assumption the number of phases per level of the QSF process is constant, say equal to $\ell+1$, for $\ell\geq 0$. Without loss of generality, we may assume the QSF process to be skip-free to the left. We will also assume that the levels are bounded to the right, since otherwise the stationary distribution cannot exist, as is shown later in this section. Without loss of generality we may then assume that they are non-negative, i.e. $m\leq0$.
Then the infinitesimal generator $\Gn$ or $q$-matrix (cf. \cite{anderson1991continuous} and \cite{kat-smit-sp-I-13}) takes the form:
\beq
\elabel{generator}
\Gn=\left[
\begin{array}{cccccccc}
\ddots  &\vdots &\vdots &\vdots &\vdots &\vdots &\vdots\\
\cdots & D & W & U^1 & U^2 & U^3 & U^{4'} \\
\cdots & 0 & D & W & U^1 & U^2 & U^{3'} \\
\cdots & 0 & 0 & D & W & U^1 & U^{2'} \\
\cdots & 0 & 0 &0 & D & W  & U^{1'} \\
\cdots & 0 & 0 & 0 &0 & D & W' \\
\end{array}
\right],
\eeq
\\
where 
the $(\ell+1) \times (\ell+1) $ sub-matrices $U^s$ $(s=1,2,\ldots)$, $W$, and $D$ represent the transition rates to the $s$-th higher level, the same level, and  the next lower level respectively. Elements of these matrices are $u^s_{ij}$, $w_{ij}$, and $d_{ij}$ where $u^s_{ij}$ is the $((m,i), (m+s,j))$ element from the matrix $U^s$, $w_{ij}$ is the $((m,i), (m,j))$ element from the matrix $W$, and $d_{ij}$ is the $((m,i), (m-1,j))$ element from the matrix $D$, for any $i,j\in\left\lbrace 0,\ldots,\ell\right\rbrace $. The matrices $U^{k'}$ are such that the row sum of $Q$ is zero. For example these matrices can be $U^{k'}=\sum_{i=k}^{\infty} U^{i}$ and $W'=W+\sum_{i=1}^{\infty} U^{i}.$

Throughout the paper we assume that the $q$-matrix $\Gn$ is irreducible, conservative, stable, and non-explosive.
Additionally, we assume that the jump rates are allowed to be  unbounded as a function of the state and that
$X$ is the minimal process. It is convenient to denote the levels by $L_m$, so that 
$
L_m=\{(m,i)\,|\, i\in\{0,\ldots,\ell\}\}\
$
and
$\state=\cup_{m}\, L_m,$ where $-\infty\leq m \leq 0$. In view of the structure of $\Gn$ given in Eq.~\eq generator/, jumps can only take place to levels $L_k$ for $k\geq m-1$, given that the current level is $L_m$. 

Since the levels are mutually exclusive, they form a partition of the state space.
We will introduce some more notation. In accordance with 
\cite{kat-smit-sp-I-13}, for any fixed $m$ we define the sub-level set of $L_m$ to be the set of states $\Lu_{m}= \cup_{k\leq m}\, L_k$,	while the set $\Lo_m=\cup_{k\geq m}\, L_k$ is the super-level set of $L_m$.
Then clearly $\{\Lu_{m-1},\Lo_m\}$
 is a partition of $\state$ for each $m$.

Suppose that the QSF process is positive recurrent. Let $\pi$ denote the (unique) stationary distribution. We denote the $m$-level
sub-vector of $\pi$ by $\pi_m$. The stationary distributions corresponding to  $\Lu_{m-1}$ and $\Lo_m$ will be denoted by  $\utilde{\pi}_{m-1}$ and $\widetilde{\pi}_m$ respectively. Then by a standard taboo decomposition, as discussed in  \cite{kat-smit-sp-I-13}, we can express $\utilde{\pi}_{m-1}$  in terms of $\widetilde{\pi}_m$  and the expected amount of time spent in the states of $\Lu_{m-1}$.

We denote by
$
\tT[m]{(k,i)(l,j)}
$  the expected amount of time spent in state $(l,j)$ without passing through states
in the super-level set $\Lo_m$, given that the system starts in state $(k,i)$, where $(k,i), (l,j)\in \Lu_{m-1}$.

Further, denote by $\tcT[][m]$  the $\left| \Lu_{m-1}\right| \times \left| \Lu_{m-1}\right| $ matrix  with elements $\tT[m]{(k,i)(l,j)}$. Write:
\begin{equation}
\tuD=\left[\begin{array}{cccc}
0&\ldots&0&0\\
\vdots&\ddots&\vdots&\vdots\\
0&\ldots&0&0\\
0&\ldots&0&D
\end{array}\right],
\end{equation}
where $0$ stands for an $(\ell+1)\times(\ell+1)$ zero matrix and $\tuD$ has dimension $|\Lo_{m}|\times|\Lu_{m-1}|$.
Then by the skip free property to the left
\beq
\elabel{stat-1}
\utilde{\pi}_{m-1}=\widetilde{\pi}_m \, \tuD \, \null_{m}\!\!{\sf T}.
\eeq
We will next express $\tcT[][m]$ in terms of $\Gn$. The elements of the $q$-matrix $\Gn$ are denoted by $\gn{(k,i)(l,j)}$, where $(k,i),(l,j)\in\state$. Denote $\gn{(k,i)}=-\gn{(k,i)(k,i)}$, the parameter
of the exponential sojourn time in state $(k,i)$.
Further denote by $\Gn[m]$ the taboo-generator with taboo set $\Lo_{m}$, restricted to the states of $\Lu_{m-1}$, with elements
$$
\gn[m]{(k,i)(l,j)}=\left\{\begin{array}{ll}
\gn{(k,i)(l,j)},\quad &(k,i),(l,j)\in \Lu_{m-1},\\
0,\quad&\mbox{otherwise},\end{array}\right.
$$
and $\gn[m]{(k,i)}=-\gn[m]{(k,i)(k,i)}$.
The taboo is therefore imposed at the time of the first jump out of a state of the QSF process on $\Lu_{m-1}$.
By assumed irreducibility, $\Gn[m]$ is the $q$-matrix of a transient, non-conservative and minimal  Markov process on the
state space $\Lu_{m-1}$. Since the number of levels is unbounded to the left, $\Gn[m]$ is of infinite dimension.

Denoting the corresponding (minimal) transition function by $\PrT[m]{}_t=\Prt[m]{t,(k,i)(l,j)}{}$ where all elements ${(k,i),(l,j)\in \Lu_{m-1}}$. It follows that 
\beq
\elabel{time}
\tT[m]{(k,i)(l,j)}=\int_0^\infty \Prt[m]{t,(k,i)(l,j)}{}dt<\infty.
\eeq
We can state the following lemma.
\begin{lem}
\label{inverse}
The inverse matrix $\Gn[m]^{-1}$ of $\Gn[m]$ exists and $\Gn[m]^{-1}=-\tcT[][m]$. Further,
$$\tcT[][m]\!\Gn[m]=\Gn[m]\,\, \tcT[][m]=-\ident_m,$$ where $\ident_m$ is the identity matrix on $\Lu_{m-1}$.\\
\end{lem}
\proof
By virtue of Anderson \cite{anderson1991continuous} (Theorem 2.2.2) the minimal transition function $\PrT[m]{}_t$ satisfies the Kolmogorov forward and backward equations
and it is the unique solution on $\state$. 

Define a scalar $a>0$. Then  $\null{}_m\!R(a)=\int_0^\infty e^{-a t}\PrT[m]{}_t \, dt$ is the associated resolvent with elements $\null_{m}\!r_{(k,i) (l,j)}(a)$, where 
$(k,i), (l,j)\in \Lu_{m-1}$.
By virtue of Anderson \cite{anderson1991continuous}, Propositions (2.1.1) and (2.1.2), it holds that
$$
a  \,  \null{}_m\!R(a)=\ident_m+\Gn[m] \, \null{}_m\!R(a)=\ident_m+\null{}_m\!R(a)\Gn[m].
$$

Hence, by virtue of the backward equation, for $(k,i), (l,j)\in \Lu_{m-1}$ it holds that
\beq
\elabel{res-2}
\begin{aligned}
a \, \null{}_m\!r_{(k,i) (l,j)}(a)&=\delta_{(k,i) (l,j)}+\gn[m]{(k,i)(k,i)} \, \null{}_m\! r_{(k,i) (l,j)}(a)\\
&+\sum_{(k',i')\ne (k,i)}\gn[m]{(k,i)(k',i')} \, \null{}_m\! r_{(k',i') (l,j)}(a),
\end{aligned}
\eeq
where Kronecker delta $
\delta_{(k,i)(l,j)}=\left\{\begin{array}{ll}
1,\quad &(k,i)=(l,j)\\
0,\quad&\mbox{otherwise}.\end{array}\right.
$

Observe that $\null{}_m\!r_{(k) (l,j)}(a)\uparrow \tT[m]{(k,i)(l,j)}$ for $a \downarrow 0$.
Using that the summation in Eq. \eq res-2/ contains only non-negative terms, 
we can now take the limit $a\downarrow 0$ in Eq. \eq res-2/ and use the monotone convergence theorem to obtain
\beq
\elabel{resolvent}
\delta_{(k,i) (l,j)}+\gn[m]{(k,i)(k,i)}\tT[m]{(k,i)(l,j)}+
\sum_{(k',i')\ne (k,i)}\gn[m]{(k,i)(k',i')}\tT[m]{(k',i')(l,j)}=0.
\eeq
Eq. \eq{resolvent}/ is equivalent to
$$
\ident_{m_,(k,i) (l,j)}+(\Gn[m]\,\tcT[][m])_{(k,i) (l,j)}=0.
$$
This yields that $\Gn[m]\, \, \tcT[][m]=-\ident_m$. 

The relation $\tcT[][m]\Gn[m]=-\ident_m$ is analogously proved, by using the forward equation for the
resolvent.\qed

As a consequence of the result above we obtain:
\beq
\elabel{stat-2}
\utilde{\pi}_{m-1}=-\widetilde{\pi}_m \tuD\Gn[m]^{-1}.
\eeq
Define $T$ to be the $(\ell+1)\times (\ell+1)$ subblock of $(-\Gn[m])^{-1}$ corresponding to the states in $L_{m-1}$. Note that without further restrictions on $D$ it is unfortunately not possible to calculate it independently of the rates of states in $\Lu_{m-1}.$
From Eq. \eq stat-1/ it follows that
$$
\pi_{m-1}=\pi_m B,\quad  B=DT,
$$
where $B_{ij}$ represents the expected local times spent in $(m-1,j)$, before absorption into $\Lo_m$, given that the process starts in $(m,i)$.

By homogeneity, we can recursively show that
\beq
\elabel{stat-3}
\pi_m=\pi_0 B^{-m}, \quad \mbox{for} \quad m\leq 0.
\eeq
Now  the existence of the stationary distribution implies that $$\sum_{m\leq 0} B^{-m}{\bf 1}_{\ell+1}=(I-B)^{-1}{\bf 1}_{\ell+1}<\infty,$$ with $I$ the $(\ell+1)\times (\ell+1)$ identity matrix and  ${\bf 1}_{\ell+1}$ the $\ell+1$-dimensional vector consisting of ones. 

The next lemma shows that it is not possible to have an unbounded state space to the right and a stationary distribution for a homogeneous QSF process.

\begin{lem}
Consider a homogeneous QSF process where the number of levels is unbounded in the negative direction. In order for the process to be positive recurrent, the number of levels has to be bounded in the positive direction.
\end{lem}
\llabel{lem:posbounded}
\proof
Suppose that the levels are not bounded to the right.
Were the stationary distribution to exist, the same reasoning as in the above (Eq.~(\ref{stat-3})) would yield that
$\pi_m=\pi_n B^{-m+n}$, $m\leq n$. Analogously, by homogeneity it would follow that $\pi_{n}=\pi_{(n,0)} {\bf a}$, with ${\bf a}$ an $(\ell+1)$-dimensional positive column vector,
independent of $n$ (see also the proof of Theorem~\ref{stat} below). 
Then 
$$
\sum_{m\leq n, j} \pi_{(m,j)}=\pi_{(n,0)}{\bf a}^T(I-B)^{-1}{\bf 1}_{\ell+1}.
$$
Taking the limit $n\to\infty$, on the one hand yields that $\sum_{m\leq n, j} \pi_{(m,j)}\to 1$. On the other hand, the right-hand side of the above equation
is the product of two scalars, $\pi_{(n,0)}$ and  ${\bf a}^T(I-B)^{-1}{\bf 1}_{\ell+1}$. 
Since $\lim_{n\to\infty}\pi_{(n,0)}=0$, the product equals 0 as well, in the limit $n\to\infty$. This is a contradiction.\qed

\subsection{Special choice of \bf{\emph{D}}}
Imposing extra structure on $D$ implies a result on the structure of the stationary distribution, that is described below. 
This extra structure covers the case when $D=c\cdot r$, where $c$ is a column vector and $r$ is a row vector, and $\cdot$ denotes matrix product (of potentially non-square matrices). In other words, the rows of $D$ are proportional to  each other. 
Let the $i^{th}$ column of $T$ be denoted by the vector $t_i=\left(
t_{0i},t_{1i},\ldots,t_{\ell i} \right)^T$, for $i\in \left\lbrace 0,\ldots,\ell \right\rbrace  $.

Then $B$ has the following form: 
$$B=DT=\left[
\begin{array}{cccc}
c_0r \inp t_0&c_0r\inp t_1&\cdots&c_0r\inp t_\ell  \\
c_1r\inp t_0&c_1r\inp t_1&\cdots&c_1r\inp t_\ell  \\
\vdots&\vdots&\ddots&\vdots\\
c_\ell r\inp t_0&c_\ell r\inp t_1&\cdots&c_\ell r\inp t_\ell  \\
\end{array}
\right]=c \cdot \hat{t},$$
where 
 $\hat{t}=(r\inp t_0,r\inp t_1,\ldots,r\inp t_\ell)$ and 
 $$
 B^n=(c\inp \hat{t})^n=c \inp \big((\hat{t}\inp c)^{n-1}\hat{t}\big)=
 \left( \sum_{i=0}^\ell c_ir\inp t_i\right)^{n-1} \!\!\!\! c\inp \hat{t},\quad  n\geq 1.
 $$
Notice that matrix $B$ has only one nonzero eigenvalue $\gamma= \hat{t}\,\inp\, c$ with left eigenvector $\hat{t}$ of multiplicity $1,$ provided that $ \hat{t}\,\inp\, c$ is finite.
Combination with Eq. \eq stat-3/ implies that the stationary distribution of the levels has a product form. 
The following theorem holds.

\begin{thm}\label{thm:colrow}
Assume that $D$ has a column times row structure as is  described above.\\
If  $\gamma<1$ and $\ell<\infty$, then the Markov process $X$ is positive recurrent.\\
If the Markov process $X$ is positive recurrent, then $\gamma<1$. In either case:
\llabel{stat}
\beq
\elabel{qsf}
\pi_{(m,j)}=\gamma^{-(m+1)} \sum_{i=0}^\ell \pi_{(0,i)}c_i \, \hat{t}_j.
\eeq
\end{thm}
\proof
Suppose that $X$ is positive recurrent. 
By  Eq.~\eq stat-3/, the stationary distribution necessarily has the form of Eq.~\eq qsf/.
Taking the summation over the levels $m$ in Eq.~\eq qsf/ yields a finite expression. This implies that
$\sum_{m\leq M}\gamma^{-(m+1)}<\infty$ for any level $M$. Hence, $\gamma<1$ necessarily. 

By separating state $(0,0)$ from level 0, one can compute the expected sojourn time in the states of $\state\setminus\{(0,0)\}$ by inverting the (negative of the) finite rate matrix restricted to this set, analogously to the above derivations.
Define $x_{(0,i)}$ to be equal to an unknown constant $x_{(0,0)}$ times the expected sojourn time in state $(0,i)$ given a start in state  $(0,0)$ before absorption into $(0,0)$, for all $i>0$.  Then define $x_{(m,j)}$ from Eq. \eq qsf/ by substituting $\pi_{(0,i)}$ by $x_{(0,i)}$ in the right-hand side.
Note that $x_{(0,0)}$ is the only unknown constant. 
Since $\gamma<1$, they can be normalised to yield a (unique) probability distribution on $\state$.
\qed

It  is still hard to compute $B$ directly, therefore we introduce even more structure on $D$ in the next section.
 More precisely, we assume that $c$ has only one non-zero component. Without loss of generality, we may assume that $c_{0}$ is the only non-zero element of $c$. 
 
\section{Exit states and successive lumpability}
\llabel{s3}

Our main assumption concerns the presence of {\em exit states}, which are defined below.
\begin{defn} Let $M\subset\state$. Then $(m, i)\in M$ is an exit state for $M$, if
\begin{description}
\item[i)]
$\sum_{(k',j')\not\in M}\gn{(k,j)(k',j')}=0$ for all $(k, j)\in M$ with $(k,j)\ne (m,i)$.
\item[ii)] $\sum_{(k', j')\not\in M}\gn{(m,i)(k',j')}>0$.
\end{description}
\end{defn}
For the remainder of this section we make the following assumption,   
which  implies that for all $m$, level $m-1$ can only be reached directly from level $m$ through state $(m,0)$.
 In other words, $(m,0)$ is an exit state to level $m-1$.

\begin{assumption}
\llabel{L3.2}
\elabel{ass1}
State $(m,0)$ is an exit state for superset $\Lo_m$, for all $m=\ldots,-1,0.$
\end{assumption}

\subsection{Stationary distribution}
\llabel{subs3}

Under assumption \ref{L3.2}, the sub-matrix $D$ contains only one non-zero row. 
This implies that matrix $B$ has only one non-zero row (the first) as well. Then by virtue of Theorem \ref {stat}:
\beq
\elabel{one-step}
\pi_{(m-1,j)}=\pi_{(m,0)} b_{0j}, 
\eeq
where $b_{0j}=c_0\,r\inp t_j$, for $j=0,\ldots,\ell$.

Let $\beta_j=b_{0j}/b_{00}$. We know: $\gamma =c_0 \hat{t}_0=c_0\, r\inp t_0=b_{00}<1$ and thus:
\beq
\elabel{productform}
\pi_{(m,j)}=\gamma^{-m}\beta_j\pi_{(0,0)},
\eeq
for $j=0,\ldots,\ell$.

As has been mentioned at the end of section \ref{s2},  explicit computation of the matrix $B$ is in general hard, because it can be of infinite dimension and it can have a complicated structure.
In the presence of an exit state however, 
it turns out to be simpler  to express 
$\pi_m$ in terms of $\utilde{\pi}_{m-1}$, cf. Eq. (\ref{rev1}).

By tabooing on the set $\Lu_{m-1}$, there exists an $|\Lu_{m-1}|\times (\ell+1)$ matrix $R_{m}$,  such that 
\beq
\elabel{rev1}
\pi_m=\utilde{\pi}_{m-1}R_m.
\eeq
By our homogeneity assumptions, $R_m$ is in fact independent of $m$, and so we suppress the dependance on $m$ in our further notation, whenever possible.
We will next define an  embedded $q$-matrix $\check{Q}$ on $L_m$ and an $\Lu_{m-1}\times L_m$ matrix $\check{A}$ by  removing the entrance state $(m+1,0)$, such that the mean (local) sojourn times in the state of $L_m$ do not change, given a start in $\Lu_m$. 
We next state the following theorem using this shift.

\begin{thm}
Consider a level homogenous process with an exit state in each level. Then:
$$
\pi_m=\utilde{\pi}_{m-1}R,
$$
where
$$
R:=-\check{A}\check{Q}^{-1},
$$
where $\check{A}$ and $\check{Q}$ are $|\Lu_{m-1}|\times (\ell+1)$  and $(\ell+1)\times(\ell+1)$ matrices with elements
\begin{equation}\label{eq:esl1}
\check{q}_{ij}=w_{ij}+\sum_{s=1}^\infty \sum_{r=0}^{\ell} u^s_{ir}\, \dfrac{d_{0j}}{\sum_{v=0}^\ell d_{0v}}, 
 \end{equation}
\begin{equation}\label{eq:esl2}
    \check{a}_{(k,i)(m,j)}=u_{ij}^{m-k}+\sum_{s=m+1-k}^\infty \sum_{r=0}^{\ell} u^s_{ir}\, \dfrac{d_{0j}}{\sum_{v=0}^\ell d_{0v}}, \mbox{ for }k< m, 0\leq i,j \leq \ell,
 \end{equation}
respectively.
\end{thm}

Notice that
$d_{0j}/\sum_{r=0}^\ell d_{0r}$ is the probability that  level set $L_{m}$ is entered at state $(m,j)$, given that a downward transition (to set $L_{m}$) occurs starting in state $(m+1,0)$.

\begin{proof}
We will explicitly compute $R$, containing  the expected sojourn times spent in the states of level $L_m$, before absorption into
$\Lu_{m-1}$, given that the process starts in $L_m$.

To this end, note that since $(m+1,0)$ is an exit state for level $\Lo_{m+1}$, it is an entrance state for the set $\utilde{L}'_{m}$, defined below:
$$
\utilde{L}'_{m:}=\Lu_m\cup\{(m+1,0)\}.
$$
This implies that $\utilde{L}'_{m}$ can only be reached from states in $\state\setminus \utilde{L}'_{m}$ through $(m+1,0)$.
The technique (successive lumping) developed in Katehakis and Smit \cite{katehakis2012successive} can be used to compute the expected sojourn  time spent in
the extended level $L'_m=L_m\cup\{(m+1,0)\}$ before absorption into $\Lu_{m-1}$, given that the process starts at this extended level $L'_m$.

We use this technique to compute  $\tilde{Q}$, the transition rate matrix of size $(\ell+2)\times(\ell+2)$ embedded on  level $L'_{m}$ for any $m$. Then $\tilde{Q}$ has elements:
\beq
\elabel{tildeQ}
\tilde{q}_{(k,i)(l,j)}=\left\{\begin{array}{ll}
w_{ij},\quad &(k,i),(l,j)\in L_{m},\\
\sum_{s=1}^\infty \sum_{r=0}^\ell u^s_{ir},\quad&  (k,i)\in L_{m}, (l,j)=(m+1,0),\\
d_{0j},\quad& (k,i)=(m+1,0), (l,j)\in L_{m},\\
-\sum_{j=0}^\ell d_{0j},\quad& (k,i)=(l,j)=(m+1,0).
\end{array}\right.
\eeq
Note that
the transitions leading to states in superlevel $\Lo_{m+1}\backslash \{(m+1,0)\}$ are mapped to the entrance state $(m+1,0)$.
By virtue of Lemma \ref{inverse} the expected sojourn time spent in  level $L'_{m}$ before absorption into the sub-level set $\Lu_{m-1}$ is obtained by inverting $-\tilde{Q}$. 
Katehakis, Smit and Spieksma \cite{kat-smit-sp-I-13}, Eq.~(9), show that
\beq 
\elabel{pi-tilde}
\pi'_m=-\utilde{\pi}_{m-1}\, \tilde{A}\tilde{Q}^{-1},
\eeq
where $\pi'_m=(\pi_m,\pi_{(m+1,0)})$ and the elements of $\tilde{A}$ are given below with $k<m$:
\beq
\elabel{tildeA}
\tilde{a}_{(k,i)(l,j)}=\left\{\begin{array}{ll}
u^{m-k}_{ij},\quad & (l,j)=(m,j)\,\, \mbox{and}\,\, i,j\in\{0,\ldots,\ell\}\\
\sum_{s=m+1-k}^\infty\sum_{r=0}^\ell u^s_{ir},\quad &(l,j)=(m+1,0).
\end{array}\right.
\eeq
The validity of this {\sl exit state lumping}  construction of Eqs. (\ref{eq:esl1}) and (\ref{eq:esl2}) is guaranteed by  Lemma~\ref{pfff} in appendix A.
\end{proof}

\paragraph{Restriction to QBD processes}
For  QBD processes, where $U^1=U$ and $U^s=0$ for $s\geq 2$, the above  derivations  impy that:
\beq
\elabel{QBD1}
\pi_m=\pi_{m-1} R,
\eeq
where $R=-U \check{Q}^{-1}.$
If $U$ is an invertible matrix, then the result above is equivalent to
\beq
\elabel{QBD}
\pi_{m-1}=-\pi_{m} \check{Q}U^{-1}.
\eeq
In this equation we have expressed $\pi_{m-1}$ in terms of $\pi_m$, similarly to Eq.~\eq stat-3/. However, $\check{Q} U^{-1}$ is explicitly computable, provided that $U$ is invertible, 
whereas $B$ may not be invertible.

From Eqns.~\eq {productform}/ and \eq {QBD1}/ we have the following result
\beq
\elabel{eigenvalue}
(1/\gamma)\beta=\beta R.
\eeq
Particularly if $U$ is an invertible matrix, then
\beq
\elabel{eigenvalue2}
\beta (-\check{Q}) U^{-1}=\gamma \beta.
\eeq

\begin{lem}
\llabel{l4.5}
If $\ell$ is finite, then $1/\gamma$ is the maximum eigenvalue of $R$ in absolute value.
\end{lem}
\proof
Let $R$ be an irreducible non-negative matrix. Then by the Perron-Frobenius theorem (for example in Seneta \cite{seneta1980non}) 
the maximum eigenvalue is positive and real. We denote this eigenvalue by $r$. Further, 
there exist unique strictly positive left and right eigenvectors. Let $x$ denote this positive right eigenvector. 
Then for all $i= 0,\ldots, \ell$:
\beq
\elabel{PF}
rx_i=\sum_j r_{ij}x_j.
\eeq

From Eq.~\eq{eigenvalue}/ it is clear that $1/\gamma$ is an eigenvalue of $R$ with $\beta,$ (where $\beta>0$) as the left eigenvector corresponding to $1/\gamma$, i.e.:
\beq
\elabel{eigen}
(1/\gamma)\beta_j=\sum_i \beta_i r_{ij}.
\eeq
Next, we will show that $1/\gamma$ is the maximum eigenvalue of $R$, thus that $1/\gamma=r$.

From Eq.~\eq {PF}/  we get
$$
\sum_i r x_i\beta_i=\sum_j\sum_i\beta_i r_{ij}x_j.
$$
Combining with Eq. \eq {eigen}/ yields
$$
r(x\inp \beta)=\frac{1}{\gamma} (x\inp \beta).
$$
Since $x\inp \beta>0$, it follows that $r=1/\gamma$.
In other words, $1/\gamma$ is the maximum eigenvalue of $R$.
\qed

In the next section we will discuss  a single server  queueing model with Coxian arrivals and batch services. 

\section{Queueing application: analysis of the \boldmath{${\it Cox}(k)/M^Y/1$} queue}\label{s4}

We consider a queueing model in which customers arrive within a maximum of $k$ exponentially  distributed phases.
The inter-arrival distribution is a Coxian distribution of order $k$, see e.g. \cite{bhulai2006value}. 

In a Coxian distribution of order $k$, 
phase $i\in\{0,\ldots,k-1\}$ lasts an exponentially distributed amount of time with parameter $\lambda_i$, at the end of which either a new customer arrives with probability $1-q_{i}$, or a new phase starts with probability $q_{i}$, where $q_{k-1}=0$.
 Upon arrival of a new customer, a new inter-arrival distribution starts.
To avoid trivialities, we assume that $q_i>0$ for $i\leq k-2$.  Hence the inter-arrival distribution is a ${\it Cox}(k)$ distribution with parameters $(\lambda_0,\ldots,\lambda_{k-1}, q_0,\ldots,q_{k-1})$.
 
We use this naming of the phases to have a more natural state space description: here, state $(m,i)$ denotes the state of the system when there are $m$ customers in the system and the $m+1$ arriving customer has completed $i$ arrival phases. We will use the random variable $C\sim Cox(k)$ to denote the inter-arrival time.

The ${\it Cox}(k)/M^{Y}/1$ queueing system can be formulated as a QSF process $X(t)$ on the state space $\textit{\textbf{S}}=\left\lbrace L_0,L_1,\ldots\right\rbrace $, where $L_m=\left\lbrace (m,0),(m,1),...,(m,k-1)\right\rbrace $ for all $m\geq 0$. 
Service occurs according to the distribution induced by $Y$; in batches of  size $j$, with $1\le j\le b$, each with probability $p_{j}.$ We will denote the probability-generating function of $Y$ as $\phi_{Y}(\cdot).$ Note that $\phi_{Y}(x)=\sum_{j=1}^{b}p_{j}x^{j}$.

We will first analyse the general ${\it Cox}(k)/M^Y/1$-queue with finitely many phases. Then we will discuss some special subcases, where the exponential phases all have the same parameter, or the probabilities of a new inter-arrival phase are all equal upto the order. We will briefly discuss the extension to the case of a Coxian distribution of infinite order.

\paragraph{Analysis of the \boldmath${\it Cox}(k)/M^Y/1$-queue}
The transition rate matrix $Q$ for this model takes the form:
\beq
\elabel{matrixQPh}
 \Gn=\left[
\begin{array}{cccccccc}
W_0&U&0&0 &0  &\cdots \\
D'_{1}&W&U&0& 0 &\cdots \\
D'_{2}&D_{1}&W&U&0 &\cdots \\
D'_{3}&D_{2}&D_{1}&W&U &\cdots \\
\vdots   &\vdots&\vdots&\vdots &\vdots& \ddots\\
\end{array}
\right],
\eeq
with $k \times k$ sub-matrices $D_{i}=p_{i}\cdot\mu I,$
$D'_{k}=\sum_{i=k}^{y} D_{i}$, 
where $I$ is the identity of size $b\times b$ and 
$$ W_0=\left[
\begin{array}{ccccc}
-\lambda_0 &q_0\lambda_0 & \cdots&\cdots  &0 \\
0& -\lambda_1 & q_1\lambda_1& \cdots &0 \\
\vdots  &\vdots &\ddots  &\vdots \\
0 &0&\cdots &-\lambda_{k-2}&q_{k-2}\lambda_{k-2}\\
 0 &0&\cdots & 0&-\lambda_{k-1}
\end{array}
\right], \quad\quad
U=\left[
\begin{array}{cccc}
 (1-q_0)\lambda_0&\cdots  &0&0 \\
(1-q_1)\lambda_1&\ddots &0&0 \\
\vdots&\ddots &\vdots&\vdots \\
(1-q_{k-2})\lambda_{k-2}& \cdots &0&0 \\
\lambda_{k-1}&  \cdots&0 &0\\
\end{array}
\right].
$$
and $W=W_0+\mu I$.
Note that $U$ has a $c\cdot r$ structure,
with
$$
c=\left(\begin{matrix} (1-q_0)\lambda_0\\ \vdots\\(1-q_{k-2}) \lambda_{k-2}\\ \lambda_{k-1}\end{matrix}\right)\quad\mbox{and } r=(1,0,\cdots,0),
$$
and therefore fits the analysis of Section~\ref{s2}. Note further  that the structure of Eq.~(\ref{matrixQPh}) is the transpose of the matrix introduced in Eq.~(\ref{generator}).  We have to interchange the roles of the matrices $U$ and $D$ in the formul\ae\ of the previous sections to fit the model under consideration. Another possible way to look at the problem is to use the original negative numbering of the levels. Then we do not need the interchanging described above, but the numbering is not induced by the model anymore.

We assume that the QSF process is ergodic. We shall prove below that for ergodicity it is necessary and sufficient to require that the mean number of arrivals per unit time is smaller than the mean number of service completions, which is specified in the relation below:

\beq  
\elabel{ERGO}
\frac{1}{\sum_{i=0}^{k-1}\prod_{l=0}^{i}q_l\lambda^{-1}_i}=\frac{1}{\Exp C}<{\mu}{\sum_{j=1}^b jp_j}=\mu \Exp Y.
\eeq
When ergodicity holds, we can use the results of Theorem \ref{thm:colrow}, and  the stationary distribution takes the form: (again, note the change of notation, since we are considering positive levels now)
$$
\pi_{(m,j)}=\gamma^{m-1}\sum_{i=0}^{k-1}\pi_{(0,i)}c_i\hat{t}_j=\gamma^{m-1}\sum_{i=0}^{k-1}\pi_{(0,i)}\lambda_i(1-q_i)\hat{t}_j,
$$
where for notational convenience, we use $q_{k-1}=0$.
The factor $\gamma$ and the
vector $\hat{t}$ denote the distribution of the phase and are implicitly given as largest positive eigenvalue and corresponding left eigenvector of the rate matrix $R$ for a fixed level. Note that once $\gamma$ and $\hat{t}$ are known, the distribution of level 0 can be deduced. 

An alternative procedure is done by observing that
the state  $(m,0)$ is an entrance state for $\Lo_m$ (and an exit state for the shifted partition $\Lo_{m-1}\cup\{(m,0)\}$).
This allows us to use the results from Section~\ref{s3} for computing the rate matrix in the reverse (downward) direction.
However, we prefer not to use the shifted partition, because the emerging
 stationary distribution in both phases and levels will then have a notationally less amenable form, as we shall see. Therefore we will stay with the current level partitions.  

In consideration of the statements above, the following relation holds (cf.~Eq.~\eq pi-tilde/) for levels $m\geq 1$:  
\beq
\elabel{iterative}
\pi_m=\pi_{m+1}R=-\widetilde{\pi}_{m+1}\, \tilde{A}\tilde{Q}^{-1},
\eeq
where $\tilde{A}$ and $\tilde{Q}$ are given in \eq tildeA/ and \eq tildeQ/. We refer to \cite{kat-smit-sp-I-13}, where a similar formula is also provided, but not the results regarding the product form.
For this model, the matrix $\tilde{Q}$ takes the form:
\begin{equation}
\elabel{Q}
\tilde{Q}=\left[
\begin{array}{cccccc}
-\lambda_0 &q_0\lambda_0&\cdots &0&0&0\\
\mu&-(\lambda_1+\mu)&q_1\lambda_1&\cdots&0&0\\
\vdots&\vdots &\ddots &\ddots&\vdots&\vdots\\
\mu&0 &\cdots& &-(\lambda_{k-2}+\mu)&q_{k-2}\lambda_{k-2} \\
\mu&0&\cdots&&0&-(\lambda_{k-1}+\mu)\\
\end{array}
\right].
\end{equation}
Note that $\tilde{Q}$ is independent of the batch size distribution.
The matrix $\tilde{A}$ has dimension $k\times bk$, and can therefore be written as:
$$
\tilde{A}=(\tilde{A}_1\ \tilde{A}_2\,\ldots \tilde{A}_b)^{T},
$$
where $\tilde{A}_i$ is the following $k\times k$ matrix: 
\begin{equation}
\elabel{A}
\tilde{A}_i=\left[
\begin{array}{cccc}
\displaystyle\sum_{j\ge i}p_j\mu&0&\cdots&0\\
\displaystyle\sum_{j>i}p_j\mu&p_i\mu&\cdots&0\\
\vdots&\vdots &\ddots &\vdots\\
\displaystyle\sum_{j>i}p_j\mu&0&\cdots&p_i\mu\\
\end{array}
\right],
\end{equation}
for $1<i<b$ and $\tilde{A}_b=p_b\cdot\mu I$.

Since the rate matrix $R$ expresses the stationary solution of lower levels in terms of that of higher ones, we cannot recursively compute the stationary distribution in the infinite capacity case. We still need to find $\gamma $ and $\hat{t}$ to do so.

Instead we will derive equations for $\gamma$ and the unknown (conditional) stationary probabilities of the phases.
To this end we consider the expression:
$$-\pi_{m}\tilde{Q}=\widetilde{\pi}_{m+1}\tilde{A}.$$
Because of the level product form solution (cf. Eq. \eq productform/)  we can rewrite this equation as:
\beq
\elabel{QSF-0}
-\pi_{m}\tilde{Q}=\sum_{i=1}^b\pi_{m+i}\tilde{A}_i=\pi_m\sum_{i=1}^b\gamma^i\tilde{A}_i.
\eeq
Using Eq.~\eq one-step/ (i.e.~$\pi_{(m,i)}=\hat{t}_{i}\pi_{(m,0)}$) to calculate the components of the vectors with $1\le i\le k-1$ on both sides of the above equality Eq.~\eq QSF-0/  yields:
\beq
\elabel{QSF1}
(\lambda_i+\mu)\hat{t}_i -q_{i-1}\lambda_{i-1}\cdot \hat{t}_{i-1}=\sum_{j=1}^b \gamma^j p_j\mu\cdot \hat{t}_i=\mu \phi_{Y}(\gamma)\cdot\hat{t}_{i}.
\eeq
Hence
\beq
\elabel{QSF-a}
\frac{\hat{t}_{i}}{\hat{t}_{i-1}}=\frac{q_{i-1}\lambda_{i-1}}{\lambda_i+\mu-\mu \phi_{Y}(\gamma)}.
\eeq
In other words, writing $\beta_{i}=q_{i-1}\alpha_i=\hat{t}_i/\hat{t}_{i-1}$, we get that $\hat{t}_i=\prod_{l=1}^iq_{l-1}\alpha_l\cdot\hat{t}_0$. Therefore the stationary distribution has the following form for $m\geq 1$ 
\begin{equation}
\elabel{Pformcox}
\pi_{(m,i)}=
\pi_{(1,0)}\gamma^{m-1}\prod_{l=1}^i q_{l-1}\alpha _{l},
\end{equation}
%
where $\alpha _i$  the following function of $\gamma$ derived from \eq QSF-a/:
\beq
\elabel{QSF-b}
\alpha _i=\frac{\lambda_{i-1}}{\lambda_i+\mu-\mu \phi_{Y}(\gamma)}.
\eeq
Note that $(\alpha _1,\ldots,\alpha _{k-1})$ only depend on $(q_0,\ldots,q_{k-2})$ through $\gamma$.

Using the balance equation for state $(m,0)$ ($m\geq1$) and Eq.~\eq QSF-b/ we obtain the following expression for $\gamma$ as a function of $(\alpha _1,\ldots,\alpha _{k-1})$:
\beq
\elabel{QSF-c}
\gamma=\frac{\alpha _1(1-q_0)\lambda_0+\alpha _1\displaystyle\sum_{i=1}^{k-1}\prod_{l=1}^{i}\alpha _lq_{l-1}(1-q_i)\lambda_i}{(\lambda_0-\lambda_1)\alpha _1+\lambda_0},
\eeq
where we have substituted $\phi_Y(\gamma)\mu=\lambda_1+\mu-\lambda_0/\alpha _1$. 
Equations \eq QSF-b/ and \eq QSF-c/ provide a system of $k$ equations in the unknowns $\gamma$ and $\alpha _i$, $i=1,\ldots,k-1$. 
As is easily checked, $\alpha _i=\lambda_{i-1}/\lambda_i$, $i=1,\ldots,k-1$, $\gamma=1$ form a solution of this system. 

Using the same balance equation, we obtain a fixpoint   equation for $\gamma$:
\beq
\elabel{QSF-c1}
\gamma=F(\gamma),
\eeq
with $F:\Re\to\Re$ given by
\beq
\elabel{QSF-c2}
F(\gamma)=\frac{\lambda_0(1-q_0)+\sum_{i=1}^{k-1}\prod_{l=1}^iq_{l-1}\frac{\lambda_{l-1}}{\lambda_l+\mu-\mu\phi_Y(\gamma)}(1-q_i)\lambda_i+\gamma\mu \phi_Y(\gamma)}{\lambda_0+\mu}.
\eeq

Again, $\gamma=1$ is a fixpoint  of $F$.
The function $F$ is tediously but  easily checked to be a convex function of $\gamma$ on the interval $[0,1+\epsilon]$ for some $\epsilon>0$. It is positive on this interval, and hence larger than the left hand side with $\gamma=0,$ i.e.~the value of $F(0)$. Computing the derivative at $\gamma=1$
yields, that this is strictly larger than 1 if and only if condition \eq ERGO/ holds. Hence, $F$ has one fixpoint  $\gamma<1$  if and only
if \eq ERGO/ holds, showing that \eq ERGO/ is  necessary and sufficient for ergodicity. This is a standard argument used for deriving ergodicity conditions by means of a probability generating function approach, see e.g. \cite{ross2013applied} Section 4.5. 

Hence solutions $\gamma$ and $\{\alpha _i\}_{i=1}^{k}$ can be determined (i) by determining the unique fixpoint  $\gamma<1$ of $F$  and then (ii) insert $\gamma$
into Eq. \ref{QSF-b}.

Next we will express the steady state probabilities of the first level in terms of the steady state probability in state $(0,0)$. For sake of presentation, we omitted the  derivation and present the results immediately:
\beq
\elabel{pi_10}
\pi_{(1,0)}=\pi_{(0,0)}\dfrac{\lambda_0}{\mu} \dfrac{1-\gamma}{1-\phi_Y(\gamma)},
\eeq
\beq
\elabel{pi_0i}
\pi_{(0,i)}=\pi_{(0,0)}\dfrac{\lambda_0}{\lambda_i}\prod_{s=0}^{i-1}q_s \left(1+\dfrac{(1-\gamma)\phi_Y(\gamma)}{\gamma(1-\phi_Y(\gamma))}\sum_{j=1}^i\prod_{l=1}^{j}\alpha _l\right),\,\,\,\, \mbox{for}\,\, i\geq 1,
\eeq
and
\beq
\elabel{pi0}
\pi_{(0,0)}=\left( \left( 1+\frac{\lambda_0}{\mu(1-\phi_Y(\gamma))}\right) \left( 1+\sum_{i=1}^{k-1}\prod_{j=0}^{i-1}\frac{\lambda_0}{\lambda_i}q_j\right)\right)^{-1}.
\eeq

Then by using Eqns.~\eq Pformcox/, \eq pi_10/, \eq pi_0i/, and \eq pi0/, for $m\geq 1$ we have 
\begin{equation}
\elabel{Pformcox0}
\pi_{(m,i)}=\pi_{(0,0)} \dfrac{\lambda_0(1-\gamma)\gamma^{m-1}}{\mu(1-\phi_Y(\gamma))}\prod_{l=1}^i q_{l-1}\alpha _{l}.
\end{equation}

We summarise all our findings above in the next theorem.
\begin{thm}
The ${\it Cox}(k)/M^{Y}/1$ queue is ergodic in and only if Eq.~\eq ERGO/ is satisfied.
If this is the case, the stationary distribution of the ${\it Cox}(k)/M^{Y}/1$ queue on levels $\geq1$ is given by  \eq Pformcox0/. The factors $\alpha _i$, $i=1,\ldots,k-1$ and $\gamma$ can be calculated from Eq.~\eq QSF-c1/ and Eq.~\eq QSF-b/. The boundary level can by found by Eq.~\eq pi_0i/ and \eq pi0/. 
\end{thm}

\begin{remark}[Finite capacity queues]\elabel{remark1}
In case of a finite capacity queue of size $S$,  the stationary distribution can be computed recursively, in terms of $\pi_S$.
By using the fact that level $S$ has an entrance state $(S,0)$ from the right, the analysis described in the previous section can be used and yields that
$$
\pi_{(S,i)}=-\tilde{Q}^{-1}_{S}\pi_{(S,0)},
$$
where now
$$
\tilde{Q}_S=\left[
\begin{array}{cccccc}
-\lambda_0&\lambda_0&&\cdots&0&0\\
\mu&-(\lambda_1+\mu)&\lambda_1&\cdots&0&0\\
\vdots&\vdots &\ddots &&\vdots&\vdots\\
\mu&0&&\cdots&-(\lambda_{k-2}+\mu)&\lambda_{k-2}\\
\mu&0&\cdots&&0&-\mu\\
\end{array}
\right].
$$
Then using \eq iterative/ with the  matrices $\tilde{Q}_m$ and $\tilde{A}_{m}$ associated with level $m<S$, one may compute the lower level stationary probabilities (up to a constant). Final renormalisation yields the correct stationary distribution.
\end{remark}
\begin{remark}[Variable service rates]
Suppose that service rates $\mu$ and batch probabilities $p_j$ are equal to $\mu^m$ and $p_j^m$ for level $m\leq S,\, j=1,...,b$, respectively, for some $S<\infty$.Then the stationary distribution of ${\it Cox}(k)/M^Y/1$ queue on the levels $S+k$, $k\geq 1$, can be computed in exactly the same manner as in the above, yielding the stationary probabilities of these levels upto a constant.
Meanwhile, the stationary distribution for levels $m\leq S$ can be computed from $\tilde{\pi}_S$ by using \eq iterative/ with the  matrices $\tilde{Q}_m$ and $\tilde{A}_{m}$ associated with level $m$, as described in Remark \ref{remark1}. Again, final renormalisation yields the correct stationary distribution.

\end{remark}

\subsection{A \boldmath${\it Cox}(k)$ inter-arrival distribution with homogeneous parameters}
In this section we consider several subcases for the model described in the previous section.
\paragraph{Homogeneous rates}
 First let us assume that the rates in the exponential distribution are all equal: $\lambda_i=\lambda$, $i=0,\ldots,k-1$. 
Then from Eq.~\eq QSF-b/ we get
\beq 
\elabel{QSF-b-hom}
\alpha _i=\frac{\lambda}{\lambda+\mu-\mu \phi_Y(\gamma)}=:\alpha ,\quad i=1,\ldots,k-1,
\eeq
in other words, the phase factors $\alpha _i$ have become independent of the phase.
Then Eq.~\eq QSF-c/ reduces to
\beq
\elabel{QSF-c3}
\gamma=\alpha (1-q_0)+\sum_{i=1}^{k-1}\big(\prod_{l=0}^{i-1}q_l\big)\alpha ^{i+1}(1-q_i).
\eeq
The function $F$ can be simplified, but we can also directly insert Eq.~\eq QSF-b-hom/ for $\alpha $ in Eq.~\eq QSF-c3/ to obtain that the solution $\gamma$
is a fixpoint  of the equation $\gamma=F^h(\gamma)$ with
$$
F^h(\gamma)=\frac{\lambda(1-q_0)}{\lambda+\mu-\mu\cdot\phi_Y(\gamma)}+\sum_{i=1}^{k-1}\Big(\prod_{l=0}^{i-1}q_l\Big)\Big(\frac{\lambda}{\lambda+\mu-\mu\cdot\phi_Y(\gamma)}\Big)^{i+1}(1-q_i)
.
$$
In the same manner as for the function $F$, we can deduce that the ergodicity condition \eq ERGO/ is necessary and sufficient for $F^h$ to have a fixpoint  $\gamma<1$. 

\paragraph{Homogeneous rates and probabilistic phase transitions}
Next we assume that additionally the probabilities of a new phase are all equal, that is, $q_i=q$, $i=0,\ldots,k-2$. 
This has no further impact on $\alpha $. However, Eq.~\eq QSF-c/ reduces further to
\beq
\elabel{QSF-c4}
\gamma=\sum_{i=0}^{k-2}q^i\alpha ^{i+1}(1-q) +q^{k-1}\alpha ^{k}.
\eeq
The function $F^h$ becomes
$$
F^h(\gamma)=\sum_{i=0}^{k-2}q^i\Big(\frac{\lambda}{\lambda+\mu-\mu\cdot\phi_Y(\gamma)}\Big)^{i+1}(1-q) +q^{k-1}
\Big(\frac{\lambda}{\lambda+\mu-\mu\cdot\phi_Y(\gamma)}\Big)^{k}.
$$
Further, the stationary distribution of the levels $m\geq1$ has a product form
\beq
\pi_{(m,i)}=\pi_{(0,0)}\dfrac{\lambda(1-\gamma)\gamma^{m-1}(q\alpha )^i}{\mu(1-\phi_Y(\gamma))},\quad m\geq 1,
\eeq  
and 
\beq
\elabel{pi00}
\pi_{(0,i)}=\pi_{(0,0)}q^i\left(1+\dfrac{(1-\gamma)\phi_Y(\gamma)}{\gamma(1-\phi_Y(\gamma))}\sum_{j=1}^i \alpha ^j\right),\quad \mbox{for}\,\, i=1,\cdots, k-1,
\eeq
where
$\pi_{(0,0)}=\dfrac{(1-q)\mu(1-\phi_Y(\gamma))}{(1-q^k)(\lambda+\mu(1-\phi_Y(\gamma)))}=\dfrac{(1-q)(1-\alpha )}{1-q^k}$.

\paragraph{Homogeneous rates and deterministic phase transitions: the $E_k/M^Y/1$-queue}
Taking $q_i=q=1$, $i=0,\ldots,k-2$, finally yields an Erlang $(k,\lambda)$ inter-arrival distribution.
Then Eq.~\eq QSF-c/ takes the very simple form
\beq
\elabel{QSF-c5}
\gamma=\alpha ^{k}.
\eeq
The function $F^h$ is  simply given by
$$
F^h(\gamma)=\Big(\frac{\lambda}{\lambda+\mu-\mu\cdot\phi_Y(\gamma)}\Big)^k.
$$
Finally, 
 the stationary distribution has a product form (except for level 0) determined by only  one factor
 \beq
\pi_{(m,i)}=\pi_{(0,0)}\dfrac{\lambda (1-\gamma)\alpha  ^{i+k(m-1))}}{\mu(1-\phi_Y(\gamma))},\quad m\geq 1,
\eeq  
and 
\beq
\elabel{pi000}
\pi_{(0,i)}=\pi_{(0,0)}\left(1+\dfrac{(1-\gamma)\phi_Y(\gamma)}{\gamma(1-\phi_Y(\gamma))}\sum_{j=1}^i \alpha ^j\right),\quad \mbox{for}\,\, i=1,\cdots, k-1,
\eeq
where
$\pi_{(0,0)}=\dfrac{1-\alpha }{k}$.

Notice that this is not surprising. The process is then essentially a one-dimensional process with a GI/M/1-structure.
This is known to have a product form stationary distribution (cf. Asmussen \cite{asmussen2003applied}). 

We further note that the ${\it Cox}(k)/M^Y/1$-queue with constant rates and probabilistic phase transitions is {\em not} of the (pure) GI/M/1-type. Indeed, the probabilities are only constant with respect to the phases with index smaller than $k-1$.

\paragraph{Distribution of the number of customers in the queue}
By taking the summation over the phases per level  we obtain  
the following modified geometric distribution for the number of customers in $E_k/M^Y/1$-queue
\beq
\elabel{pi-geom}
\bar{\pi}_m^k=\left\{\begin{array}{ll}
\displaystyle\frac{\lambda(1-\gamma)^2}{k\mu\cdot(1-\phi_Y(\gamma))} \gamma^{m-1},\quad & m>0,\\\noalign{\vspace*{4pt}}
\displaystyle 1-\frac{\lambda(1-\gamma)}{k\mu\cdot(1-\phi_Y(\gamma))},\quad &  m=0,
\end{array}\right.
\eeq
where $\bar{\pi}_m^k=\sum_{i=0}^{k-1} \pi_{(m,i)}$, $m\geq0$.

\subsection{The \boldmath${\it Cox}(\infty)/M^Y/1$-queue}

Allowing infinitely many phases still fits our framework. We get the natural extensions of formulas Eq.~\eq QSF-c/ and Eq.~\eq QSF-c2/ with $k=\infty$.
Clearly the expression for $\alpha _i$ as a function of $\gamma$ is not affected by the amount of phases.

Restricting to the case of homogeneous rates and probabilities, so that $\lambda_i=\lambda$, $q_i=q$, $i=0,\ldots$, we get the following results.

The queueing system is ergodic if and only if
\beq
\elabel{ERGO2}
\lambda \left(\frac{1-q}{q}\right)<\mu\sum_{j=1}^b jp_j.
\eeq

Then Eq.~\eq QSF-c/ becomes 
\beq
\elabel{QSF-c6}
\gamma=\frac{\alpha (1-q)}{1-q\alpha },
\eeq
and so again we obtain
$$
F^h(\gamma)=\frac{\lambda(1-q)}{\lambda(1-q)+\mu(1-\phi_Y(\gamma))}.
$$

\vspace*{0.5\baselineskip}\noindent
We summarise our results in the next theorem.
\begin{thm} 
Under the ergodicity condition \eq ERGO2/, the stationary distribution of the ${\it Cox}(\infty)/M^{Y}/1$ queue on levels $\geq1$ is given by  \eq Pformcox0/ where $\alpha $ and $\gamma$ can be calculated from Eq.~\eq QSF-b-hom/ and Eq.~\eq QSF-c6/. The boundary level can by found by Eq.~\eq pi_0i/, where
$
\elabel{pi_00}
\pi_{(0,0)}=(1-q)(1-\alpha ).$
\end{thm}

\subsection{Numerical Analysis}
Assume an ergodic ${\it Cox}(k)/M^Y/1$ queue of finite order.
In order to calculate the solutions $\{\alpha _i,i=1,\ldots,k-1;\gamma\}$, 
one can solve \eq QSF-c1/ directly, e.g. by the Newton-Raphson method, and then determine $\alpha _i$, $i=1,\ldots,k-1$, from \eq QSF-b/.

As $\gamma$ is a fixpoint  of the map $F$ associated with \eq QSF-c1/, another possibility is to approximate the fixpoint  $\gamma<1$ by selecting
$\gamma_0$ and by iteratively computing $\gamma_{n+1}=F(\gamma_n)$.
It is simply checked that the fixpoint  $\gamma=1$ is not stable, but that the fixpoint  of interest, smaller than 1, is a stable fixpoint .
Hence we can use the following scheme to approximate the desired values $\alpha _i,i=1,\ldots,k-1$, $\gamma$.

\vspace*{0.5\baselineskip}\noindent
{\bf Approximation scheme 1}
\begin{enumerate}
\item Choose $\gamma<1$; 
\item iteratively put $\gamma:=F(\gamma)$, till desired convergence;  compute $\alpha _i$ from \eq QSF-b/, $i=1,\ldots,k-1$.
\end{enumerate}

In the case of homogeneous rates, we have shown that the level factor $\gamma$ is the unique fixpoint  smaller than 1 of the function $F^h$ .
This leads to the following adapted scheme,  that we only formulate for the case of homogeneous rates.

\vspace*{0.5\baselineskip}\noindent
{\bf Approximation scheme 2 for the case $\lambda_i=\lambda$, $i\leq k-1$ (where $k=\infty$ is allowed)}
\begin{enumerate}
\item Choose $\gamma<1$;
\item iteratively put $\gamma:=F^h(\gamma)$, until desired convergence, and compute $\alpha $ from \eq QSF-b-hom/.
\end{enumerate}

It is outside the scope of the paper to discuss the rate of convergence of the scheme, as well as a detailed stopping criterion. Table 1 below shows the non-surprising property that $\gamma$ and $\alpha $ are non-increasing in $q$. Further note that even for a high value of the continuation probability $q$, $\gamma$ is already
approximately constant starting from around 20 phases.

\begin{center}
\footnotesize
\begin{tabular}{|c| c  c  |c c|c c| c c|c c|}
\hline
 $q$&$k=2$&&$k=5$ &   & $k=20$ &  & $k=1000$& &$k=\infty$ 
 &
\\
& $\gamma$ & $\alpha $ & $\gamma$ & $\alpha $ & $\gamma$ & $\alpha $ & $\gamma$ & $\alpha $ 
 &$\gamma$ & $\alpha $  \\
\hline\noalign{\vspace*{0.1cm}}\hline

0.1 & 0.4168&0.4414&0.4153 & 0.4411 & 0.4153 & 0.4411 & 0.4153 & 0.4411&0.4153&0.4411\\
\hline

0.2 & 0.3847&0.4339&0.3788 & 0.4325 & 0.3788 & 0.4325 &
0.3788 & 0.4325&0.3788&0.4325\\
\hline

0.3 & 0.3538&0.4272&0.3406 & 0.4246 & 0.3406 & 0.4246 &
0.3406 & 0.4246&0.3406&0.4246\\
\hline

0.4 & 0.3239&0.4214&0.3006 & 0.4173 & 0.3005 & 0.4172 &
0.3005 & 0.4172&0.3005&0.4172\\
\hline

0.5 & 0.2948&0.4163&0.2585 & 0.4105 & 0.2582 & 0.4104 &
0.2582 & 0.4104&0.2582&0.4104\\
\hline

0.6 & 0.2663&0.4117&0.2141 & 0.4042 & 0.2134 & 0.4042 &
0.2134 & 0.4042&0.2134&0.4042\\
\hline

0.7 & 0.2385&0.4076&0.1673 & 0.3986 & 0.1658& 0.3984 &
0.1658 & 0.3984&0.1658&0.3984\\
\hline

0.8 & 0.2113&0.4039&0.1176 & 0.3935 & 0.1147& 0.3932 &
0.1147 & 0.3932&0.1147&0.3932\\
\hline

0.9 & 0.1845&0.4006&0.0648 & 0.3890 & 0.0598 & 0.3886 &
0.0598 & 0.3886&0.0598&0.3886\\
\hline

1 & 0.1581&0.3976&0.0085 & 0.3851 & $\approx 0$ & 0.3846 &
$\approx 0$ & 0.3846&$\approx 0$&0.3846\\
\hline
\end{tabular}
\end{center}
\begin{center}
${\vtop{\hbox{Table 1. $(\lambda, \mu, \gamma, p_1, p_2, p_3)=(0.5, 0.8, 0.35, 0.25, 0.5, 0.25)$.}}}$
\end{center}

\section{Monotonicity properties and relation with the \boldmath{$D/M^Y/1$}-queue}
\llabel{s5}

\subsection{Monotonicity properties}
\llabel{subs5-1}

Next we will study monotonicity properties for the homogeneous ${\it Cox}(k)/M^Y/1$ queues, in the following sense.  If $\lambda_i=\lambda$ and $q_i=q$ for
$i=0,\ldots,k-2$, then we denote the corresponding $k$-order Coxian distribution by ${\it Cox}(k,\lambda,q)$.

In the remainder of the paper we assume that the ergodicity condition  Eq. \eq ERGO/ is satisfied.
It is the aim to 
 compare the stationary distribution of the number of customers in the system for ${\it Cox}(k,\lambda,q)/M^Y/1$ queues, with different inter-arrival distributions with the same mean inter-arrival times and associsated probabilities $q$, but with a different amount of phases.

Let $\lambda^*$ and phase probability $q$ and the number of phases $k$ be given.
Then, in order for the mean inter-arrival time to equal $1/\lambda^*$, the parameter $\lambda_k$ of the homogeneous Coxian distribution has to be equal to:
$$
\frac{1}{\lambda^*}=\frac{1}{\lambda_k}\sum_{l=0}^{k-1} q^l=\frac{1-q^{k}}{\lambda_k(1-q)},
$$
that is
$$
\lambda_k=\frac{\lambda^*(1-q^k)}{1-q}.
$$
Denote the corresponding factors in the stationary distribution (see previous section) of the corresponding QSF process by $\alpha _k$ and $\gamma_k$ respectively, and denote the stationary distribution of the number of customers in the system by $\bar{\pi}^k$.

It does not seem possible to stochastically compare these queueing systems for a different number of phases:
Table 2 below shows that there is a lack of monotonicity in the parameter $\gamma_k$, especially for high values of the continuation parameter $q$.

\begin{center}
\footnotesize
\begin{tabular}{|c| c  c  |c c|c c| c c|c c|c c |}
\hline
 $q$&$k=2$ &   & $k=5$ &  & $k=10$& &$k=20$ 
 &&$k=50$& &$k=1000$&\\
& $\gamma$ & $\alpha $ & $\gamma$ & $\alpha $ & $\gamma$ & $\alpha $ & $\gamma$ & $\alpha $ 
 &$\gamma$ & $\alpha $  &$\gamma$ & $\alpha $ \\
\hline\noalign{\vspace*{0.1cm}}\hline

0.1 &  0.4485 & 0.4734&0.4502&0.4764&0.4502&0.4764 &0.4502&0.4764&0.4502&0.4764&0.4502&0.4764\\
\hline

0.2 & 0.4439&0.4939&0.4501&0.5057&0.4502&0.5058&0.4502&0.5058&0.4502&0.5058&0.4502&0.5058\\
\hline

0.3 & 0.4371&0.5120&0.4494&0.5383&0.4502&0.5391&0.4502&0.5391&0.4502&0.5391&0.4502&0.5391\\
\hline

0.4 & 0.4286&0.5283&0.4472&0.5738&0.4502&0.5771&0.4502&0.5771&0.4502&0.5771&0.4502&0.5771\\
\hline

0.5 & 0.4189&0.5429&0.4415&0.6111&0.4499&0.6206&0.4502&0.6209&0.4502&0.6209&0.4502&0.6209\\
\hline

0.6 & 0.4082&0.5563&0.4300&0.6488&0.4482&0.6670&0.4502&0.6718&0.4502&0.6718&0.4502&0.6718\\
\hline

0.7 & 0.3968&0.5685&0.4105&0.6853&0.4413&0.7243&0.4499&0.7316&0.4502&0.7319&0.4502&0.7319\\
\hline

0.8 & 0.3848&0.5798&0.3811&0.7197&0.4195&0.7796&0.4463&0.8012&0.4502&0.8037&0.4502&0.8037\\
\hline

0.9 & 0.3725&0.5902&0.3407&0.7514&0.3658&0.8307&0.4139&0.8738&0.4484&0.8905&0.4502&0.8912\\
\hline

\end{tabular}
\end{center}
\begin{center}
${\vtop{\hbox{Table 2. $(\lambda^*, \mu, \gamma, p_1, p_2, p_3)=(0.5, 0.8, 0.35, 0.25, 0.5, 0.25)$.}}}$
\end{center}

However, in the case of deterministic number  phase transitions, an Erlang queue, there exists more structure regarding the change of $\gamma$ with respect to a increasing number of phases $k.$

\paragraph{Monotonicity properties for the \boldmath{$E_k/M^Y/1$}-queue}
Let us next restrict to the case $q=1$, in other words, the case of an Erlang inter-arrival distribution.
Then the arrival rate in the $k$-phase system is given by $\lambda_k=\lambda^* k$. Write $\rho=\lambda^*/\mu$.

The following results hold.
\begin{thm}
\llabel{gamma3}
\begin{description}
\item[a)]
The sequence $\gamma_k$ is strictly decreasing in $k$. It has a limit $\gamma^*=\lim_{k\rightarrow \infty} \gamma_k$, which is the unique    solution $\xi$ smaller than 1 of the equation below.
\beq
\elabel{fix}
\xi=e^{-(1-\phi_Y(\xi))/\rho}.
\eeq
\item[b)] The map $k\mapsto \bar{\pi}^{k}_0=1-\rho(1-\gamma_k)/(1-\phi_Y(\gamma_k))$ is a strictly decreasing function if and only if $\Prob\{Y=1\}<1$. If $\Prob\{Y=1\}=1$, i.e. the batch size equals 1 with probability 1, then $\bar{\pi}^k_0=1-\rho$, for $k=1,2,\ldots$
\end{description}
\end{thm}

\proof 
By combination of \eq QSF-b-hom/ and \eq QSF-c5/ we obtain:
\beq
\elabel{gamma2}
1/\gamma_k=\left(1+ \dfrac{(1-\phi_Y(\gamma_k))}{k\rho}\right)^{k}.
\eeq
Define 
$
g_k(x):=\left(1+ \dfrac{x}{k}\right)^{k},\quad k=0,1,\ldots.
$

Clearly this function is increasing in $x$. To show that $g_{k}$ is also increasing in $k$, we expand the expression using the binomial formula:
\begin{eqnarray}
g_{k+1}(x)&=&\sum\limits_{i=0}^{k+1}
\dbinom{k+1}{i}
\left( \dfrac{x}{\rho(k+1)}\right)^i \nonumber\\
\elabel{f1}
&=&1+\sum\limits_{i=1}^{k}\dfrac{(x/\rho)^i}{i!}\cdot\dfrac{k}{k+1}\cdots \dfrac{k+2-i}{k+1}+\dfrac{(x/\rho)^{k+1}}{(k+1)^{k+1}}.
\end{eqnarray}

On the other hand, 
\begin{eqnarray}
g_{k}(x)&=&\sum\limits_{i=0}^{k}
\dbinom{k}{i}
\left( \dfrac{x}{\rho(k)}\right)^i \nonumber\\
\elabel{f2}
&=&1+\sum\limits_{i=1}^{k}\dfrac{(x/\rho)^i}{i!}\cdot\dfrac{k-1}{k}\cdots \dfrac{k+1-i}{k}.
\end{eqnarray}

Eq.~\eq{f1}/ and Eq.~\eq{f2}/ yield via a term by term comparison that 
\beq
\elabel{f3}
g_{k+1}(x)>g_k(x), \quad \mbox{for all} \quad x>0.
\eeq
In other words, $g_{k+1}(1-\phi_Y(\gamma))>g_k(1-\phi_Y(\gamma))$, with $g_{k+1}(1-\phi_Y(1))=g_k(1-\phi_Y(1))=1$ and $g_{k+1}(1-\phi_Y(0))=(1+\rho^{-1})^{k+1}>
(1+\rho^{-1})^{k}$. Recall that $\gamma_i$ is the unique fixpoint  of the equation $1/\gamma=g_i(1-\phi_Y(\gamma))$,  with $\gamma_i\in(0,1)$, $i=k,k+1$.
For $x<\gamma_i$, $1/x>g_i(1-\phi_Y(x))$ and for $x\in(\gamma_i,1)$, necessarily $1/x<g_i(1-\phi_Y(x))$, $i=k,k+1$. Hence $\gamma_{k+1}<\gamma_k$.

Since $\gamma_k$ are non-increasing, and bounded below, the sequence has a limit, $\gamma^*$ say, with $\gamma^*<1$.
This limit solves \eq fix/ by the standard limiting argument that $\lim_{k\to\infty}(1+x/k)^k=e^x$.
The function $\gamma\mapsto e^{-(1-\phi_Y(\gamma))}$ is a convex function, with fixpoint s $\gamma=1$ and $\gamma^*<1$, derivative larger than 1 at $\gamma=1$,
and positive value at $\gamma=0$. The result then follows in a standard manner, thus completing the proof of a).

Part b) follows from the fact that 
\beq
\elabel{monoton}
\frac{1-\gamma}{1-\phi_Y(\gamma)}=\frac{1}{\sum_{i=1}^bp_i\sum_{j=0}^{i-1}\gamma^j },
\eeq
which is strictly increasing in $\gamma<1$ 
if and only if $\Prob\{Y=1\}<1$. 
\qed

Clearly, with increasing $k$, the variance of the inter-arrival time decreases. Hence, the average server utilisation strictly improves with decreasing variance, 
in the case of batch sizes $Y$, with $\Prob\{Y=1\}<1$. Whereas, if the batch size equals 1 with probability 1, the average server utilisation is equal to $\rho$ and thus constant.

We can say a little more.
We define $L_k=\sum_m m\cdot \bar{\pi}_m^k$ to be the mean number of customers in the system under the stationary distribution. Further, we denote by $W_k$ the expected sojourn time and by $V_k$ the variance. 
The following comparison result holds.
\begin{thm}
\llabel{sojourn}
The following are true:
\beq
\elabel{variance}
L_{k+1}\leq L_k \quad W_{k+1}\leq W_k, \quad \mbox{and} \quad V_{k+1}\leq V_k,
\quad \mbox{for} \quad k=1,2,\ldots.
\eeq
\end{thm}

\proof 
The expectations with respect to $\bar{\pi}_{k+1}$ and $\bar{\pi}_k$ are equal to $L_{k+1}$ and $L_k$ respectively.  

It is easy to check that
\beq
\elabel{expected}
L_k=\dfrac{\rho}{1-\phi_Y(\gamma_k)},\quad \mbox{for} \quad k=1,2,\ldots.
\eeq

Since $\gamma_k$ are non-increasing, then
$$
1-\phi_Y(\gamma_k)\leq 1-\phi_Y(\gamma_{k+1}).
$$
It follows that $L_{k+1}\leq L_k.$
Application of Little's formula yields $W_{k+1}\leq W_k$.

From \eq{pi-geom}/ it is clear that the number of customers in the $E_k/M/1$-queue is (almost) geometrically distributed. This implies that  $V_k=\rho^2\gamma_k/(1-\gamma_k)^2$. It is easy to check that this yields $V_{k+1}\leq V_k$. 

Next, we need to check the variance for a general batch size distribution:
\beq
\elabel{variance2}
V_k=\dfrac{\rho(1-\gamma_k)^2}{1-\phi_Y(\gamma_k)}\sum_{m\geq 1}m^2 \gamma_k^{m-1}-L_k^2.
\eeq
Applying geometric series sum yields
\begin{eqnarray}
\sum_{m\geq 1}m^2 \gamma_k^{m-1}&=&\gamma_k\sum_{m\geq 1}m(m-1)\gamma_k^{m-2}+\sum_{m\geq 1}m\gamma_k^{m-1}\nonumber\\
\elabel{variance3}
&=&\dfrac{2\gamma_k}{(1-\gamma_k)^3}+\dfrac{1}{(1-\gamma_k)^2}.
\end{eqnarray}

By combination of \eq expected/, \eq variance2/, and \eq variance3/
\begin{eqnarray}
V_k&=&\dfrac{\rho(1-\gamma_k)^2}{1-\phi_Y(\gamma_k)}\left( \dfrac{2\gamma_k}{(1-\gamma_k)^3}+\dfrac{1}{(1-\gamma_k)^2}\right) -\dfrac{\rho^2}{(1-\phi_Y(\gamma_k))^2}\nonumber\\
\elabel{variance4}
&=&\dfrac{\rho}{1-\phi_Y(\gamma_k)}\left( \dfrac{1+\gamma_k}{1-\gamma_k}-\dfrac{\rho}{1-\phi_Y(\gamma_k)}\right).
\end{eqnarray}

We define $h(\gamma):=\dfrac{1+\gamma}{1-\gamma}-\dfrac{\rho}{1-\phi_Y(\gamma)}$. 

Since $
\gamma_k\geq \gamma^*$ and $ \gamma_k \downarrow \gamma^*$ by virtue of Theorem~\ref{gamma3} (a), it is sufficient to show that $h'(\gamma)=\frac{dh(\gamma)}{d\gamma}\geq 0$ for $\gamma\in[\gamma^*,1)$.
$$
h'(\gamma)=\dfrac{1}{(1-\gamma)^2}\left( 2-\rho\left( \dfrac{1-\gamma}{1-\phi_Y(\gamma)}\right)^2 \phi'_Y(\gamma)\right). $$

Applying Eq.~\eq monoton/ yields
\beq
\elabel{diff}
h'(\gamma)=\dfrac{1}{(1-\gamma)^2}\left( 2-\rho\left( \dfrac{1}{\sum_{j=1}^b p_j\sum_{i=0}^{j-1}\gamma^i}\right)^2\sum_{j=1}^b jp_j\gamma^{j-1}\right).
\eeq

Since $j\gamma^{j-1}\leq \sum_{i=0}^{j-1}\gamma^i$ for $\gamma\in [\gamma^*,1)$, it follows that
$$
h'(\gamma)\geq \dfrac{1}{(1-\gamma)^2}\left( 2- \dfrac{\rho}{\sum_{j=1}^b p_j\sum_{i=0}^{j-1}\gamma^i}\right).
$$

For $\gamma\in [\gamma^*,1)$,
$$
\frac{\rho}{\sum_{j=1}^b p_j\sum_{i=0}^{j-1}\gamma^i}\leq \frac{\rho}{\sum_{j=1}^b p_j\sum_{i=0}^{j-1}(\gamma^*)^i}.
$$
So, to show that $h'(\gamma\geq 0)$ for $\gamma\in [\gamma^*,1)$, it is sufficient to show that
\begin{eqnarray}
2\geq \frac{\rho}{\sum_{j=1}^b p_j\sum_{i=0}^{j-1}(\gamma^*)^i}= \frac{\rho}{\sum_{j=1}^b p_j\frac{1-(\gamma^*)^j}{1-\gamma^*}}, \nonumber
\end{eqnarray}
or
\beq
\elabel{gamma*}
2(1-\phi_Y(\gamma^*))/\rho\geq 1-\gamma^*.
\eeq

By virtue of Theorem~\ref{gamma3} (a), $
\gamma^*\geq 1-(1-\phi_Y(\gamma^*))/\rho$ (since $e^{-x}\geq 1-x$, for $x\geq 0$).

This implies that 
\beq
\elabel{variance5}
(1-\phi_Y(\gamma^*))/\rho\geq 1-\gamma^*.
\eeq

So that \eq gamma*/ follows.

Then
\beq
\frac{\rho}{1-\phi_Y(\gamma_k)}h(\gamma_k)\geq \frac{\rho}{1-\phi_Y(\gamma_{k+1})}h(\gamma_{k+1}).\nonumber
\eeq
In other words, we have proved that $V_k\geq V_{k+1}$.
\qed

Interestingly enough, for general batch size distribution, the stationary distribution does not stochastically decrease with increasing $k$. This follows immediately from
the fact that $\sum_{m\geq 1}\pi^k_m$ is strictly increasing in $k$, whereas $\gamma_k$ is strictly decreasing.
However, if the batch size is indentically equal to 1, then the stationary distribution has  a stochastically monotonic behaviour as a function of $k$.
\begin{cor}
\llabel{t-mon}
Suppose that $\Prob\{Y=1\}=1$.
The following is true for $k=1,2,\ldots$:
$$
\bar{\pi}_{k+1} \stless \bar{\pi}_k,
$$
or equivalently: $\sum_{m\geq M}\bar{\pi}_{m}^{k+1}\leq\sum_{m\geq M} \bar{\pi}_{m}^k$, for all $M=0,1,\ldots$.
\end{cor}

In this particular case we can also prove that $\alpha _k$ is strictly increasing in $k$.
\begin{lem}
Suppose that $\Prob\{Y=1\}=1$.
The sequence of parameters $\alpha _{k}$ is strictly 
increasing in $k$.
\end{lem}
\proof
From Eq.~\eq QSF-b-hom/ we have
$
\alpha _k=\dfrac{\lambda_k}{\lambda_k+\mu-\mu\gamma_k}.
$

We define 
\beq
f_k(x)=\dfrac{\mu}{\lambda_k}x^{k+1}-(1+\dfrac{\mu}{\lambda_k})x+1.
\eeq

Taking the first derivative of $f_k$ yields 
\beq
f'_k(x)=\dfrac{\mu(k+1)}{\lambda_k}x^{k}-(1+\dfrac{\mu}{\lambda_k}).
\eeq
 
Then $\alpha _{k}^*=\left( \frac{1+\frac{\lambda_k}{\mu}}{(k+1)}\right) ^{1/k}$ is the point where the polynomial $f_{k}(x)$ has a unique minimum.

Thus 
\beq
\elabel{beta}
\alpha _{k}\in\left( \frac{1}{1+\frac{\mu}{\lambda_k}},\alpha _{k}^*\right).
\eeq

For any $k \geq 1$, we have
$
\alpha _{k}^*=
\left( \frac{1+\frac{\lambda_k}{\mu}}{(k+1)}\right) ^{1/k}=\left(1- \frac{k}{k+1}(1-\rho)\right) ^{1/k}.
$

We define 
\beq
g(k)=
\left(1- \frac{k}{k+1}(1-\rho)\right) ^{1/k}.
\eeq

We will show that $g(k)$ is stricly increasing in $k$.
\beq
g'(k)=g(k)\left( -\dfrac{1}{k^2}\log\left(1- \frac{k}{k+1}(1-\rho)\right)- \frac{1-\rho}{k(k+1)(1+k\rho)}
 \right).
 \eeq
 
Since $g(k)=\alpha ^*_k>0$, it is sufficient to show that 
\beq
\elabel{log1}
-\dfrac{1}{k^2}\log\left(1- \frac{k}{k+1}(1-\rho)\right)- \frac{1-\rho}{k(k+1)(1+k\rho)}> 0.
\eeq

We know that:
$\dfrac{1}{k+1}>\dfrac{1}{(k+1)(1+k\rho)},
$ for $k\geq 1$.

This implies
\begin{eqnarray}
 e^{-\frac{1-\rho}{(k+1)(1+k\rho)}}>e^{-\frac{1-\rho}{k+1}}>\left(1- \frac{k}{k+1}(1-\rho)\right) ^{1/k}\nonumber\\
 \elabel{log2}
   -\dfrac{1}{k^2}\log\left(1- \frac{k}{k+1}(1-\rho)\right)> \frac{1-\rho}{k(k+1)(1+k\rho)}.
\end{eqnarray}

By combining \eq log1/ and \eq log2/ it is clear that $g(k)$ is stricly increasing in $k$. This means that $\alpha _{k}^*<\alpha _{k+1}^*$ and so $\alpha _{k+1}^*\uparrow 1$, for $k\to\infty$.

Next, for any positive $k$ and $x\in (0,1)$, the functions $f_k$ and $f_{k+1}$  have $2$ intersection points: at $0$ and $1$. This follows from the relation:
\beq
\elabel{difference}
f_{k+1}(x)-f_k(x)=\dfrac{1}{k(k+1)\rho}x(1-x)^2\sum_{i=1}^{k}i x^{i-1}.
\eeq

This also implies that $f_{k+1}(x)>f_k(x)$ for all $x\in(0,1)$. Hence $\alpha _{k}<\alpha _{k+1}<1$.

In other words we can say that $\alpha _{k}$ is strictly increasing in $k$ and the proof is complete.
\qed

\subsection{Comparison of batch service queues with  Erlang arrivals versus a deterministic inter-arrival time}
\llabel{subss5-2}
Let us assume for the moment that the batch size is identically equal to 1, i.e. $\Prob\{Y=1\}$.
Consider the $D/M/1$-queue with mean inter-arrival time $1/\lambda^*$ and mean service time $1/\mu$, where $\lambda^*=\lambda_k/k$.
In this case it is quite well-known that the stationary distribution of the $E_k/M/1$-queue, with the same mean inter-arrival and mean service time distribution, converges setwise and weakly to the
stationary distribution of the $D/M/1$-queue, as $k\to\infty$.
We will generalise these results when there are batch services.
By virtue of  Asmussen \cite{asmussen2003applied} and Bhat \cite{bhat2008introduction} the stationary distribution $\bar{\pi}_D$ of the $D/M/1$-queue is given by
$$
\bar{\pi}_{m}^D=\left\{\begin{array}{ll}
(1-\sigma)\rho \sigma^{m-1},\quad & m>0\\
1-\rho,\quad &  m=0,
\end{array}\right.
$$
where $\sigma$ is the unique root smaller than 1 of
\beq
\elabel{eksponen}
\sigma=e^{-(1-\sigma)/\rho},
\eeq
 when $\rho<1$.
 By virtue of Theorem~\ref{fix}, this means that $\sigma=\gamma^*=\lim_{k\to\infty}\gamma_k$, and in particular
 $\sigma<\gamma_k$.
It follows directly that:
\beq
\elabel{sdm1}
\bar{\pi}^D\stless\bar{\pi}^{k+1}\stless\bar{\pi}^k,\quad k=1,2,\ldots
\eeq

Define $L_D$, $W_D$, and $V_D$ as the mean number of customers, the expected sojourn time, and the variance of customers of the $D/M/1$-queue respectively, under the stationary distribution. We now derive the result below.
\begin{cor}
\elabel{cor:DEk}
The following are true:
\begin{description}
\item[i)]
 Eq.~\eq sdm1/ holds for all $k=1,\ldots$, and
 \item[ii)]
$
L_k\downarrow L_D, \quad  W_k\downarrow W_D, \quad  \mbox{and} \quad V_k\downarrow V_D,\quad \mbox{as} \quad k\to\infty,
$
monotonically.
\end{description}
\end{cor}
Intuitively it seems clear that a similar result holds as well for the $E_k/M^Y/1$- and $D/M^Y/1$-queues, with the same mean inter-arrival times and the same  batch service distributions. So far, we have not been able to find any result on   (setwise and weak) convergence  of the stationary distribution of the $E_k/M^Y/1$-queue to the
stationary distribution of the $D/M^Y/1$-queue, although it should be completely similar to convergence results in the case of batch size equal to 1.

\paragraph{General batch size}
Suppose  now again that $\Prob\{Y=1\}<1$. 

\begin{thm}
\label{thm-D}
For the stationary distribution $\bar{\pi}^D$ of the $D/M^Y/1$-queue  the following holds:
\begin{description}
\item[i)] $\bar{\pi}_m^k\to\bar{\pi}_m^D$, for $m=0,\ldots$ and
\item[ii)]
\beq
\elabel{pi-geom1}
\bar{\pi}_{m}^D=\left\{\begin{array}{ll}
\displaystyle\frac{\rho(1-\sigma)^2}{1-\phi_Y(\sigma)} \sigma^{m-1},\quad & m>0,\\\noalign{\vspace*{4pt}}
\displaystyle 1-\frac{\rho(1-\sigma)}{1-\phi_Y(\sigma)},\quad &  m=0,
\end{array}\right.
\eeq
where $\sigma$ is the unique root smaller than 1 of 
\beq
\elabel{sigma}
\sigma=e^{-(1-\phi_Y(\sigma))/\rho}, \quad \phi_Y(\sigma)=\sum_{j=1}^b p_j\sigma^j.
\eeq
\end{description}
\end{thm}

\proof
Since the proof is quite standard, we only provide a sketch of the proof.

First, notice that the embedded processes on arrival instants are of the GI/M/1-type considered in \cite[pp. 82--86]{karl66}.
As has been derived there, it follows for the corresponding stationary distribution $\bar{\pi}^D_A$ 
of the $D/M^Y/1$-queue that
$$
\bar{\pi}^D_{A,m}=(1-\sigma)\sigma^m,
$$
with $\sigma$ the unique root smaller than 1 of Eq.~\eq sigma/.


Secondly, by using semi-regeneration, cf. Asmussen~\cite[Ch. VII.5 and p. 283]{asmussen2003applied}, we then obtain formula Eq.~\eq pi-geom1/  for the 
stationary distribution of the non-embedded $D/M^Y/1$-queue. The desired convergence properties then follow from the explicit formulae for the respective stationary distributions.
 \qed
 
It directly follows from Theorem~\ref{thm-D}(i) and Theorem~\ref{sojourn}
that the result of Corollary~\ref{cor:DEk} holds in the batch service case as well.
\begin{cor}
The monotonicity result in Corollary~\ref{cor:DEk} (ii) holds in the case of batch service.
\end{cor}
As a consequence, the mean number of customers, expected sojourn time and variance are minimized by deterministic inter-arrival times. This is a well-known result
for non-batch systems (cf. Asmussen \cite[pp. 336-339]{asmussen2003applied}). Results on other performance measures are given in this section as well.

\section*{Acknowledgement}
This research was supported by Directorate General of Human Resource for Science, Technology and Higher Education of the Ministry of National Education of Indonesia.

\appendix
\section{Appendix}
\begin{lem}
\llabel{exist}
Let $X$ be a possibly  non-conservative, non-explosive, transient, stable Markov process on state space $\state$ with q-matrix $Q$. Then $Q^{-1}$ exists and $Q^{-1}=T$ where the $(i,j)$-th element of $T$ is defined as follows:
$\tau_{i,j}=\int_0^\infty p_{t,(i,j)} dt<\infty,$  where $p_{t,(i,j)}$ are the elements of the transition function $P_t$.
\end{lem}
\proof
A generalization of Lemma~\ref{inverse} and the proof is analogous.
\qed
\begin{lem}
\label{pfff}
Let $X$ be a possibly non-conservative, non-explosive, transient, stable Markov process on state space $\state$ with q-matrix $Q$. Let $s\in\state$ be given. Consider the transition rate matrix $\tilde{Q}$ on $X\backslash \{s\}$ where the elements are given by:
\begin{equation*}
\tilde{q}_{ij}=q_{ij}+q_{is}\frac{q_{sj}}{q_{s}},
\end{equation*}
where $q_s=-q_{ss}$.

Let ${q}_{ij}^{-1}$, $\tilde{q}_{ij}^{-1}$   denote the $(i,j)$-th elements of matrices $Q^{-1}$ and $\tilde{Q}^{-1}$, respectively.  
 Then the following are true for $i,j\neq s$:
\begin{itemize}
\item $\tilde{q}_{ij}^{-1}=q^{-1}_{ij},$
\item $q_{sj}^{-1}=\sum_{r\ne s} \frac{q_{sr}}{q_{s}}\tilde{q}_{rj}^{-1}.$
\end{itemize}
\end{lem}
\proof
By Lemma~\ref{exist} the inverse matrices $Q^{-1}$ and $\tilde{Q}^{-1}$ exist, and the entries are equal to the expected sojourn time spent in each state of $\state$ and $\state_s=\state\backslash\{s\}$ respectively.
It is convenient to use a representation based on  
the  jump chain on $\state$ with transition matrix $P^J$. The entries of $P^J$ are given by:
\beq
\elabel{jump1}
p_{i,j}^J=\frac{q_{i,j}}{q_{i}},
\eeq
where $q_{i}=-q_{i,i}$, for $i,j\in \state.$

Clearly $P^J$ is a sub-stochastic matrix, since $Q$ is transient. We denote its $n$-th iterate by $P^{J,n}$. The $0$-th iterate is the identity.
It follows from Anderson \cite{anderson1991continuous}, Proposition (4.1.1) (see also Eq.~\eq time/ of this paper) that:
\beq
\elabel{res-3}
\tau_{i,j}=\sum_{n\geq0} p_{i,j}^{J,n}\frac{1}{q_{j}}.
\eeq
The jump transition probability matrix $\tilde{P}^J$ associated with $\tilde{Q}$ on $\state_s$, is given by\beq
\elabel{jump2}
\tilde{p}_{i,j}^J=p_{i,j}^J+p_{i,s}^J p_{s,j}^J.
\eeq
This is precisely the transition matrix of the jump chain associated with $Q$, \emph{embedded on} $\state_s$.
By using  Eq.~\eq {res-3}/ and the fact that the $(i,j)$-th element of $-Q^{-1}$ represents the expected amount of time spent in state $j$, given a start in state $i$ before absorption outside $S$, we directly obtain that $-\tilde{Q}^{-1}$ is the expected amount of time spent in $\state_s$ before absorption outside this set. 

This proves the first statement. For the second statement, we invoke \eq res-3/. Then for $j\ne s$
$$
\tau_{s,j}=\sum_{r\ne s}p_{s,r}^J\tau_{r,j}=\sum_{r\ne s}\frac{q_{sr}}{q_s}\tau_{r,s}.
$$
The result follows.
\qed

\bibliographystyle{plain}
\bibliography{References}

\end{document}